\documentclass[12pt]{article}
\usepackage{color}
\usepackage{amsmath}
\usepackage{amsmath, amssymb, amscd, amsthm, amsfonts}
\usepackage{extarrows}
\usepackage{geometry}
\usepackage{authblk} 
\usepackage{hyperref} 

\usepackage[utf8]{inputenc}
\newenvironment{reproof}{{\noindent\it Proof of Theorem \ref{polynomial}.}}{\hfill $\square$\par}
\usepackage{enumerate}
\usepackage{enumitem}
\usepackage{verbatim}
\newtheorem{theorem}{Theorem}
\newtheorem{proposition}[theorem]{Proposition}
\newtheorem{lemma}[theorem]{Lemma}

\newtheorem{corollary}[theorem]{Corollary}

\newtheorem{definition}[theorem]{Definition}

\newcounter{cases}
\newcounter{subcases}[cases]
\newenvironment{mycase}
{
    \setcounter{cases}{0}
    \setcounter{subcases}{0}
    \newcommand{\case}
    {
        \par\indent\stepcounter{cases}\textbf{Case \thecases.}
    }
    
}
{
    \par
}
\renewcommand*\thecases{\arabic{cases}}

\newcommand{\ratio}{\frac{\omega(T_\omega)}{\omega(G)}}
\newcommand{\ratioF}{\frac{\omega(T_\omega)}{\omega(F)}}

\newcommand{\ratioohat}{\frac{\omega^{\ast}(T_{\omega^{\ast}})}{\omega^\ast (G)}}

\title{The Minimum Weighting Ratio Problem and Its Application in Chordal Graphs}

\author{
    {\small\bf Hui Lei$^1$}\thanks{email:  hlei@nankai.edu.cn}\quad
    {\small\bf Mei Lu$^2$}\thanks{email: lumei@tsinghua.edu.cn}\quad
     {\small\bf Yongtang Shi$^3$}\thanks{email:  shi@nankai.edu.cn}\quad
     {\small\bf Jian Sun$^3$}\thanks{email:  \text{jian\_sun@nankai.edu.cn}}\quad
    {\small\bf Xiamiao Zhao$^2$}\thanks{Corresponding author: email:  zxm23@mails.tsinghua.edu.cn}\\
    {\small 1. School of Statistics and Data Science, LPMC and KLMDASR, Nankai University, Tianjin, China.}\\
    {\small 2. Department of Mathematical Sciences, Tsinghua University, Beijing 100084, China.}\\
    {\small 3. Center for Combinatorics and LPMC, Nankai University, Tianjin, China.}
}
\date{}
\begin{document}
\maketitle


\begin{abstract}
Constructing the maximum spanning tree $T$ of an edge-weighted connected graph $G$ is one of the important research topics in computer science and optimization, and the related research results have played an active role in practical applications. In this paper, we are concerned with the ratio of the weighted sum of a spanning tree $T$ of $G$ to the weighted sum of $G$, which we try to minimize. We propose an interesting theorem to simplify this problem and show that this optimal problem can be solved in polynomial time. Furthermore, we apply the optimal problem in chordal graphs.
\end{abstract}

\textbf{Keywords:} maximum spanning tree; ratio; polynomial time; chordal graph

\section{Introduction}
In an edge-weight undirected graph $G(V,E,\omega)$ with vertex set $V$, edge set $E$ and $\omega$ is an edge-weight function $\omega: E(G)\to \mathbb{R}_{\geq 0}$. The minimum (maximum) spanning tree problem aims to find a tree $T$ which spans $V$ with minimum (maximum) weights $\omega(T):=\sum_{e\in T}\omega(e)$. Spanning tree is one of the main research topics in discrete mathematics, which has important theoretical and practical significance.

Typically, researchers are more concerned with the least weighted structure that connects all the vertices, i.e., minimum spanning tree. It has been proved that the minimum spanning tree problem can be solved in polynomial time by using Prim's Algorithm \cite{prim1957shortest} and Kruskal's Algorithm \cite{kruskal1956shortest}. The follow-up research mainly focuses on designing faster algorithms and finding spanning trees with special structures or properties \cite{DBLP:journals/jacm/Chazelle00a,DBLP:journals/combinatorica/GabowGST86,DBLP:journals/siamcomp/KatohIM81,karger1995randomized,DBLP:journals/tcs/LiFJZ18}. More specifically, March et al. considered the Euclidean minimum spanning tree problem and proposed a fast algorithm with theoretical analysis and applications \cite{DBLP:conf/kdd/MarchRG10}. Elkin studied the problem of constructing a minimum spanning tree in a distributed network and he devised a protocol that constructs the minimum spanning tree in $\widetilde{\Omega}(\mu(G,\,\omega)+\sqrt n)$ rounds, where $\mu(G,\,\omega)$ is the minimum spanning tree radius of the graph \cite{DBLP:journals/jcss/Elkin06}. Compared to the minimum spanning tree, there are relatively few research results on the maximum spanning tree \cite{DBLP:journals/comgeo/AgarwalMS91,DBLP:journals/jal/Gavril87,DBLP:journals/tcs/GalbiatiMM97}. In this paper, we mainly focus on the maximum spanning tree problem. More specifically, we are interested in the weighting ratio of the maximum spanning tree and the total graph, denoted as $t_\omega(G):= \omega(T_\omega)/\omega(G)$, where $T_\omega$ represents the maximum spanning tree of $G$ with edge-weight function $\omega$, and $\omega(G) = \sum_{e\in E(G)}\omega(e)$.

 Logistics transportation plays an important role in industrial production today, and constructing an economical and efficient transportation network is of great concern. We construct a graph with vertices representing cities, an edge connects two vertices if some mode of transportation can directly reach the two represented cities, and the weight on the edge is the corresponding cost. In this way, the graph model can represent the relevant transportation network design problems. More specifically, the demand to build an economical transport network can be seen as assigning weights to the edges (the cost of transportation modes between the two cities) in the graph to minimize the ratio of the weight of the spanning tree to the weight of the whole graph. 

The parameter $t_\omega(G)$ of $G$ proposed above describes the proportion of the weight required to connect all vertices to the total weight in the worst case. Thus it measures the cost of maintaining the connectivity of a graph. Then we can focus on the graph problem to minimize the ratio $t_\omega(G)$.

When all edges are weighted $1$, i.e., the edge-weight function is $\omega^1 : E(G) \to \{0,1\}$,  we define the uniform ratio $u(G)$ of $G$ as
$$u(G):= \frac{\omega^1(T_{\omega^1})}{\omega^1(G)} = \frac{|V(G)|-1}{|E(G)|}.$$
With respect to the weighting ratio $t_{w}(G)$, we obtain a quite interesting result as follows:
\begin{theorem}\label{gapth}
For any $\epsilon > 0 $, there exists a graph G and $\omega: E(G)\to \mathbb{R}_{\geq 0}$, such that $$\frac{t_\omega(G)}{u(G)} <\epsilon.$$
\end{theorem}
{\it Proof}    Let $T$ be a tree with $|V(T)| = t+1, t \geq 1$.
    Consider a graph $G$ consisting of a complete graph $K_\kappa$ and  $T$ joined by a vertex. 

    We set all edges in $K_\kappa$ with weight $1$, and edges in $T$ with weight $0$. Thus we have that:
    \begin{equation*}
        \frac{t_\omega(G)}{u(G)} = \frac{(\kappa-1)/{\kappa\choose 2}}{(\kappa+t-1)/({\kappa\choose 2}+t)}=\frac{2}{\kappa} \cdot \frac{\frac{\kappa^2-\kappa}{2}+t}{\kappa-1+t}.
    \end{equation*}
    Let $t=\kappa^2$. When $\kappa$ is large enough, we have
    \begin{equation*}
    \frac{2}{\kappa} \cdot \frac{\frac{\kappa^2-\kappa}{2}+t}{\kappa-1+t} \leq \frac{2}{\kappa}\cdot 1 = \frac{2}{\kappa}.
    \end{equation*}
    So when $\kappa > \frac{1}{2\epsilon}$, $\frac{t_\omega(G)}{u(G)}< \epsilon$ holds.  
\qed

As we have claimed above, the cost of a maximum spanning tree reflects the structural price required to connect all vertices (excluding cycles) in extreme cases. When given an edge-weight graph $G(V, E, \omega)$, it is easy to calculate $t_\omega(G)$. However, Theorem \ref{gapth} indicates that different edge-weight functions may cause a big difference to the parameter $t_\omega(G)$ even in the same connected graph.

Let $e(G)$ be the number of edges and $v(G)$ be the number of vertices. For a vertex set $S$, $G[S]$ is denoted as the induced graph of $G$ by $S$. In this paper, we aim to find the edge-weight function that minimized $t_\omega(G)$. We call this problem \textbf{minimum weighting ratio problem}. 
Let $W$ denote the set of non-negative edge-weight functions on $E(G)$. The minimum weighting ratio can be written as  $\inf_{\omega \in W} t_\omega(G)$. With deep analysis, we prove that for any graph with $n$ vertices, the minimum weighting ratio problem can be solved in polynomial time on $n$.
\begin{theorem}\label{polynomial}
    The problem $\inf_{\omega \in W}t_\omega(G)$ can be solved in polynomial time.
\end{theorem}
We can also obtain that the minimum value of $t_\omega(G)$ is closely related to some local properties of $G$. 
\begin{theorem}\label{inthm}
    Let $G$ be a graph. We have
    $\inf_{\omega\in W} t_\omega(G) = \min_{S\in Co(G)}\frac{|S|-1}{e(G[S])}$, where $Co(G):=\{S\subset V(G)~|~ G[S] \text{ is connected }\}.$
\end{theorem}
For graphs with specific properties, we can characterize the structure of the optimal solution.

The chordal graph is a special class of graphs with numerous applications such as semidefinite optimization and evolutionary trees \cite{enright2011application}, which belongs to the class of perfect graphs. An undirected graph $G$ is chordal if, for any cycle with length of at least four, there exists an edge connecting two nonadjacent vertices in the cycle. Equivalently, let $C_\kappa$ denote the cycle of size $\kappa$, then $G$ is a chordal graph if $G$ is induced $C_\kappa$-free for any $\kappa\geq 4$. With many interesting properties, chordal graph is one of the most studied graph classes in graph theory and has extensive applications in various fields \cite{DBLP:journals/tcs/BaiX24,DBLP:journals/jct/HurlbertK11,DBLP:conf/iwoca/PapadopoulosT22}. Numerous algorithms have been developed to explore how to find a chordal subgraph with large size \cite{balas1986fast,bhowmick2015new}.

For a given graph with $n$ vertices and $m$ edges, there are many related results in studying its maximum chord subgraph. In the 1980s, Erd\"{o}s, Gy\'{a}rf\'{a}s, Ordman, and Zalcstein \cite{erdos1989size}
 proved that when $m= n^2/4+1$, the size of maximum chordal subgraph is at least $3n/2-1$. In 2022, Gishboliner and Sudakov  \cite{gishboliner2023maximal} found that any graph with $n$ vertices and ${n^2(\kappa-1)}/{2\kappa}+a$ edges has a chordal subgraph with at least $(\kappa-1/\kappa)n-O(\sqrt{n})$ edges, where $1\leq a\leq n^2/{(2\kappa^2+2\kappa)}$.

 In this paper, we give the optimal solution to the minimum weighting ratio problem in chordal graphs and state that the minimum ratio can be bounded by the clique number, i.e., the size of the maximum clique. Since then, the number of edges of a chordal graph can be bounded by clique number through the inequality $t_\omega(G)\leq \frac{v(G)-1}{e(G)}$. We can give a lower bound about the clique number of maximum chordal subgraphs for any fixed graph.

In Section \ref{section-complexity}, we will prove Theorem \ref{polynomial}. In Section \ref{section_chordal graph}, we characterize the structure of the optimal solution to the minimum weighting ratio problem in the chordal graph.

\section{Computational Complexity}\label{section-complexity}
In this section, we further specify the optimal solution based on its structural characteristics and prove that the problem can be solved in polynomial time through in-depth analysis. Lemma \ref{the2} plays an important role in the proof.

\begin{lemma}\label{lemma1}
    Assume that $$\frac{x_1}{y_1}\leq\frac{x_2}{y_2}\leq\frac{x_3}{y_3}\leq \cdots \leq \frac{x_n}{y_n},$$
    where $x_i \geq 0$, and $y_i >0$. Then
    \begin{equation}\label{ieq1}
        \frac{x_1}{y_1} \leq \frac{x_1+x_2+\cdots+x_n}{y_1+y_2+\cdots+y_n} \leq \frac{x_n}{y_n}.
    \end{equation}
\end{lemma}
{\it Proof}  
    We just prove the left part of the inequality. By induction on $n$. When $n=2$,
    \begin{equation*}
        \frac{x_1}{y_1} - \frac{x_1+x_2}{y_1+y_2}=\frac{x_1y_2-x_2y_1}{y_1(y_1+y_2)}\leq 0
    \end{equation*}
    holds since $\frac{x_1}{y_1}\leq \frac{x_2}{y_2}.$
    Suppose that inequality (\ref{ieq1}) holds when $n=k-1$. Then we have
    \begin{equation*}
        \frac{x_1}{y_1}\leq\frac{x_2}{y_2}\leq\frac{x_2+x_3+\cdots+x_n}{y_2+y_3+\cdots+y_n}.
    \end{equation*}
    Then by the inequality on $n=2$, we have that
    \begin{equation*}
        \frac{x_1}{y_1}\leq \frac{x_1+x_2+\cdots+x_n}{y_1+y_2+\cdots+y_n}. \quad \qed
    \end{equation*}

\begin{lemma}\label{le2}
    Let $a_1$, $a_2$, $b_1$ and $b_2$ be four positive real numbers, $f(x) = \frac{a_1+b_1x}{a_2+b_2x}$. If $$\frac{b_1}{b_2}\geq f(x)$$ for all $x\geq 0$, then $f'(x)\geq 0 $.
\end{lemma}
{\it Proof}  
This can be proved easily by taking the derivative of $f(x)$.  
\qed

Next, we consider the problem of minimizing $t_\omega(G)$ for a fixed graph $G$. Initially, the value region of the edge-weight function $w$ is $\mathbb{R}_{\geq 0}$, but the next Lemma shows that we can reduce the region to $\{0,1\}$. Let $W^{0-1}=\{\omega|\omega: E(G) \to \{0,1\}\}$ and $W=\{\omega|\omega: E(G) \to \mathbb{R}_{\geq 0}\}$. The key Lemma of this section is as follows.
\begin{lemma}[0-1 Weighting Function Lemma]\label{the2}
Let $G$ be a graph. Then
\begin{equation*}
    \inf_{\omega \in W}t_\omega(G) = \min_{\omega_{0-1}\in W^{0-1}}t_{\omega_{0-1}}(G).
\end{equation*}
\end{lemma}
{\it Proof}  
    Let $\omega^{\ast}$ be the optimal solution of $\inf_{w\in W}t_{w}(G)$ with the smallest number of different edge weights. Denote $\{a_1,\ldots,a_q\}$ be the set of corresponding edge weights. We may assume that $0\leq a_1< a_2<\cdots<a_q\leq 1 $ since we can divide $a_i$ $(i=1,\ldots,q)$ by the maximum value of $\omega^{\ast}$ such that $t_{\omega^{\ast}}(G)$ remains unchanged. Recall $T_{\omega^\ast}$ represents the maximum spanning tree of $G$ with edge-weight function $\omega^\ast$.
    For $1\leq i\leq 2$, let

    $A_i^0=\{e\in E(G)\setminus E(T_{\omega^\ast}):\omega^\ast(e)=a_i\}$,

    $A_i^1=\{e\in E(T_{\omega^\ast}):\omega^\ast(e)=a_i\}$,

    $A_i=\{e\in E(G):\omega^\ast(e)=a_i\}$, that is $A_i=A_i^1\cup A_i^0$.

    We first show that $q \leq 2$. Suppose $ q\geq 3$. Then there exists  $i$ ($2\leq i \leq q-1$) such that $a_{i-1}<a_i<a_{i+1}$. Let $b$ be a real number in $[a_{i-1},a_{i+1}]$.
    Define $\omega_b: E(G) \to \mathbb{R}_{\geq 0}$ with $\omega_b(e)=\omega^{\ast}(e)$ (resp. $\omega_b(e)=b$) if $\omega^{\ast}(e)\not=a_i$ (resp. $\omega^{\ast}(e)=a_i$) for any $e\in E(G)$. Then $ \omega_b\in W$. We will show that $T_{\omega_{b}}=G[E(T_{\omega^{\ast}})]$ is a maximum spanning tree under the edge-weight function $\omega_b$, where $G[E(T_{\omega^{\ast}})]$ is a subgraph of $G$ induced by $E(T_{\omega^{\ast}})$.


%

Let $T_k$ and $T^{\prime}_k$ be the forests in the $k$-th step of Kruskal algorithm with $\omega^{\ast}$ and $\omega_b$ as edge-weight function, respectively. It suffices to prove that $T_k = T^{\prime}_k$, for $k=1,2,\ldots n$.  Since $a_i<a_q$ and $b<a_q$, we have $T_1 = T^{\prime}_1$. Assume that $T_k = T^{\prime}_k$.  Let $e_{k+1}$ be the edge added to $T_k$ in $(k+1)$-th step, with $\omega^{\ast}(e_{k+1}) = a_j$.
\begin{mycase}
    \case\label{case1} $a_j \geq a_{i+1}$ or $a_j \leq a_{i-1}$, say $a_j \geq a_{i+1}$. In this case, if $\omega^{\ast}(e)>a_j$ for any $e\in E(G)\setminus E(T_k)$, then $T_k+e$ is not a forest.
     Then $e_{k+1}$ is one of the largest weighted edge that can be put in $T^{\prime}_k$ which implies  $T^{\prime}_{k+1}= T_{k+1}$.
    \case $a_j = a_i$.
    In this case,  $e_{k+1}$ is one of the edge with largest weight we can put into $T_{k}$ by the same argument as in Case \ref{case1}. Since $\omega^{\prime}(e) = b\geq a_{i-1}$, $e$ is also the largest-weighted edge that can be added into $T^{\prime}_k$. So $T_{k+1}=T^{\prime}_{k+1}$.
\end{mycase}

Therefore $T_{\omega^{\ast}} = T_{\omega_{b}}$ holds. Note that
\begin{equation*}
    \frac{\omega^{\ast}(T_{\omega^{\ast}})}{\omega^{\ast}(G)} = \frac{|A_1^1|a_1+\cdots+|A_i^1|a_i+\cdots+|A_q^1|a_q}{|A_1|a_1+\cdots+|A_i|a_i+\cdots+|A_q|a_q}.
\end{equation*}
If $\frac{|A_i^1|}{|A_i|}\geq \frac{\omega^{\ast}(T_{\omega^{\ast}})}{\omega^{\ast}(G)},$  by letting $b=a_{i-1}$, we have
\begin{equation*}
\begin{aligned}
    \frac{\omega_{a_{i-1}}(T_{\omega_{a_{i-1}}})}{\omega_{a_{i-1}}(G)} &= \frac{|A_1^1|a_1+\cdots+|A_i^1|a_{i-1}+\cdots+|A_q^1|a_q}{|A_1|a_1+\cdots+|A_i|a_{i-1}+\cdots+|A_q|a_q}\\
    &\leq \frac{|A_1^1|a_1+\cdots+|A_i^1|a_i+\cdots+|A_q^1|a_q}{|A_1|a_1+\cdots+|A_i|a_i+\cdots+|A_q|a_q}\\
    &= \frac{\omega^{\ast}(T_{\omega^{\ast}})}{\omega^{\ast}(G)}.
\end{aligned}
\end{equation*}
The inequality holds from Lemma \ref{le2}. So 
$\omega_{a_{i-1}}$ is also an edge-weight function that minimizes the weighting ratio but the  number of different edge weights is smaller than $\omega^{\ast}$,  a contradiction with the choice of $\omega^{\ast}$.

Similarly, if $\frac{|A_i^1|}{|A_i|}\leq\frac{\omega^{\ast}(T_{\omega^{\ast}})}{\omega^{\ast}(G)},$ then
we can let $b=a_{i+1}$ and derive a contradiction as above.
So $q\leq2$.

Assume $q=2$. Then $\frac{\omega^{\ast}(T_{\omega^{\ast}})}{\omega^{\ast}(G)}=\frac{|A_1^1|a_1+|A_2^1|a_2}{|A_1|a_1+|A_2|a_2}.$
If
$\frac{|A_1^1|}{|A_1|}\geq \frac{\omega^{\ast}(T_{\omega^{\ast}})}{\omega^{\ast}(G)},$  then we let $\omega': E(G) \to \{0,1\}$ with $\omega'(e)=0$ (resp. $\omega'(e)=1/a_2$) if $\omega^{\ast}(e)=a_1$ (resp. $\omega^{\ast}(e)=a_2$) for any $e\in E(G)$. So
$$\frac{\omega'(T_{\omega'})}{\omega'(G)}=\frac{|A_2^1|}{|A_2|}\leq\frac{\omega^{\ast}(T_{\omega^{\ast}})}{\omega^{\ast}(G)},$$ where $E(T_{\omega'})=E(\omega^{\ast})$ and we are done.
If
$\frac{|A_1^1|}{|A_1|}<\frac{\omega^{\ast}(T_{\omega^{\ast}})}{\omega^{\ast}(G)},$ then we let $\omega': E(G) \to \{0,1\}$ with $\omega'(e)=1$  for any $e\in E(G)$. So
$$\frac{\omega'(T_{\omega'})}{\omega'(G)}=\frac{|A_2^1|}{|A_2|}\leq\frac{\omega^{\ast}(T_{\omega^{\ast}})}{\omega^{\ast}(G)},$$ where $E(T_{\omega'})=E(\omega^{\ast})$ and we are done.

If $q=1$, then we let $\omega': E(G) \to \{0,1\}$ with $\omega'(e)=1$  for any $e\in E(G)$. Then $\frac{\omega'(T_{\omega'})}{\omega'(G)}=\frac{\omega^{\ast}(T_{\omega^{\ast}})}{\omega^{\ast}(G)},$ and we are done. 
\qed

 For a given edge-weight function $\omega\in W^{0-1}$, let 
 $supp(\omega)=\{e\in E(G)~|~\omega(e)=1\}$. Without loss of generality, we assume that the maximum spanning tree is constructed by Kruskal's Algorithm in the following paragraphs. Then we can obtain the following two lemmas.
 \begin{lemma}\label{le4}
   Let $\omega_{i}\in W^{0-1}$ and $T_{\omega_i}$ be the maximum spanning tree of edge-weight function $\omega_i$, where $i=1,2$. If $supp(\omega_2)=supp(\omega_1)\cup\{e\}$ with $e\notin supp(\omega_1)$, then
   $$\omega_{2}(T_{\omega_2})-\omega_{1}(T_{\omega_1})=\left\{
   \begin{array}{ll}
   0 & \mbox{if $e\notin E(T_{\omega_{2}})$},\\
   1 &  \mbox{if $e\in E(T_{\omega_{2}})$}.
   \end{array}\right.$$
 \end{lemma}
{\it Proof}   Obviously, $\omega_{2}(T_{\omega_2})-\omega_{1}(T_{\omega_1})\ge 0$. Let $e\notin supp(\omega_1)$.
 Note that when constructing the maximum spanning tree via Kruskal algorithm, unless only $e$ satisfies the selection criterion, we are more likely to choose the other edges that satisfy the condition.

 If $e\notin E(T_{\omega_{2}})$, then $\omega_{2}(T_{\omega_2})=\omega_{1}(T_{\omega_2})\leq \omega_{1}(T_{\omega_1})$. On the other hand, $\omega_{2}(T_{\omega_2})\geq \omega_{2}(T_{\omega_1})=\omega_{1}(T_{\omega_1})$. Hence $\omega_{2}(T_{\omega_2})=\omega_{1}(T_{\omega_1})$ and we are done.

 Now we assume $e\in E(T_{\omega_{2}})$. Assume that $T_{\omega_{1}}$ is constructed in the order of edge selection $$e_1,e_2,\ldots,e_j,e_{j+1},\ldots,e_{n-1},$$
in which $e_{j+1}$ is the first edge with weight of $0$ during the construction of $T_{\omega_{1}}$. Then $e\notin\{e_1,\ldots,e_j\}$. Since $supp(\omega_2)=supp(\omega_1)\cup\{e\}$ and $e\in E(T_{\omega_{2}})$, $T_{\omega_{2}}$ can be assumed being constructed in the order of edge selection $$e_1,e_2,\ldots,e_j,e,e_{j+2}^{\prime},\ldots,e_{n-1}^{\prime}.$$
According to the edge selection criterion in the construction of maximum spanning tree stipulated above, 
it can be derived that $\omega_{2}(e_{j+k}^{\prime})=0$ when $2\leq k\leq n-1-j$. Otherwise, $F=G[\{e_1,\ldots,e_j,e_{j+k}^{\prime}\}])$ is a forest with
$\omega_{1}(F)=\omega_{2}(F)>\omega_{1}(T_{\omega_1})$, a contradiction. Thus we have $\omega_{2}(T_{\omega_2})-\omega_{1}(T_{\omega_1})=1$. 
\qed

\begin{lemma}\label{le5}
    Assume $supp(\omega_1)\subseteq supp(\omega_2)$. Then for any $ e\in E(G)\setminus supp(\omega_2)$, we have
    $$\omega_{1}^{\prime}(T_{\omega_{1}^{\prime}})-\omega_{1}(T_{\omega_1})\geq \omega_{2}^{\prime}(T_{\omega_{2}^{\prime}})-\omega_{2}(T_{\omega_2}),$$
    where $supp(\omega_{i}^{\prime})=supp(\omega_{i})\cup\{e\}~ (i=1,2)$.
\end{lemma}
{\it Proof}   Let $\Delta_1=\omega_{1}^{\prime}(T_{\omega_{1}^{\prime}})-\omega_{1}(T_{\omega_1})$ and $\Delta_2=\omega_{2}^{\prime}(T_{\omega_{2}^{\prime}})-\omega_{2}(T_{\omega_2})$.
  If $e\notin E(T_{\omega_{2}^{\prime}})$, then $\Delta_2=0$ and the result holds by Lemma \ref{le4}.

  If $e\in E(T_{\omega_2^{\prime}})$, then $\Delta_2=1$ by Lemma \ref{le4}. We may assume that the orders of edge selection in the process of constructing  $T_{\omega_{1}^{\prime}}$ and $T_{\omega_{2}^{\prime}}$ are
  $$e_1,e_2,\ldots,e_{j-1},e_j,e_{j+1},\ldots,e_{n-1},$$
  and
  $$e_1,e_2,\ldots,e_{j-1},e_{j}^{\prime},e_{j+1}^{\prime},\ldots,e_{n-1}^{\prime}$$
  respectively, and $e_j$ is the first edge that appears in $T_{\omega_1^\prime}$ but not in $T_{\omega_2^\prime}$. 
  The difference that $T_{\omega_{1}^{\prime}}$ selects edge $e_{j}$ while ${\omega_{2}^{\prime}}$ selects edge $e_{j}^{\prime}$ in $j$-th step implies that
  $$\omega_{2}^{\prime}(e_{j}^{\prime})>\omega_{2}^{\prime}(e_{j}),$$
  which means $\omega_{2}^{\prime}(e_{j}^{\prime})=1$ and $\omega_{2}^{\prime}(e_{j})=\omega_{1}^{\prime}(e_{j})=0$. Thus for all $0\leq k\leq n-1-j$, we have  $\omega_{1}^{\prime}(e_{j+k})=0.$

  If $e\in \{e_{j}^{\prime},e_{j+1}^{\prime},\ldots,e_{n-1}^{\prime}\}$, then  $F=G[\{e_1,e_2,\ldots,e_{j-1},e\}]$ is a forest in $G$. Thus  $\omega_{1}^{\prime}(F)>\omega_{1}^{\prime}(T_{\omega_{1}^{\prime}})$, a contradiction. Therefore $e\in\{e_1,e_2,\ldots,e_{j-1},e_{j}\}$  which implies  $e\in T_{\omega_{1}^{\prime}}$. Thus we have that $\Delta_{1}=1$ by Lemma \ref{le4}.

  Then for the both subcases, we prove that $\Delta_{1}\geq \Delta_{2}.$ 
\qed
Based on the result of Lemma \ref{le5}, we observe that $\max_{T\in\mathcal{T}}\omega(T)$ is a submodular function of $\omega$ over $W^{0-1}$. Then we can prove the computational complexity of the minimum weighting ratio problem is polynomial.

\begin{reproof}
Based on Lemma \ref{the2}, we only need to prove that the problem $$\min_{\omega^{0-1}\in W^{0-1}}t_{\omega^{0-1}}(G)$$ can be solved in polynomial time.

By Lemmas \ref{le4} and  \ref{le5}, we conclude that the problem $\min_{\omega^{0-1}\in W^{0-1}}t_{\omega^{0-1}}(G)$ is exactly to minimize the ratio of a
submodular function and modular function. Thus it can be solved in polynomial time
by a combinatorial algorithm \cite{fleischer2003push}.
\end{reproof}

\section{The Minimum Weighting Ratio of the Chordal Graph}\label{section_chordal graph}
The $0$-$1$ Weighting Function Lemma indicates that when minimizing the ratio $t_\omega(G)={\omega(T_\omega)}/{\omega(G)}$, we can directly determine which edges weight $1$, then remaining edges weight $0$. 
Now we can prove Theorem \ref{inthm}. For convenience, we rewrite it here.
\begin{theorem}
    Let $G$ be a graph. We have
    $\inf_{\omega\in W} t_\omega(G) = \min_{S\in Co(G)}\frac{|S|-1}{e(G[S])}$, where $Co(G):=\{S\subset V(G)~|~ G[S] \text{ is connected }\}.$
\end{theorem}
{\it Proof}   Let $\omega\in W$.
By Lemma \ref{the2}, we may assume that $\omega\in W^{0-1}$. Let $F=\{e\in E(G)~|~\omega(e)=1\}$ and $H=G[F]$. Then we have $$t_\omega(G)=\frac{\omega(T_\omega)}{\omega(G)}= \frac{v(H)-C(H)}{e(H)},$$
where $C(H)$ denotes the number of components of $H$. 

Let $S_i$ be the vertex set of components in $H$, where $1\le i\le C(H)$. Since $e(H)\leq \sum_{i=1}^{C(H)} e(G[S_i])$, we have $$\frac{v(H)-C(H)}{e(H)}\geq\frac{\sum_{i=1}^{C(H)}(|S_i|-1)}{\sum_{i=1}^{C(H)}e(G[S_i])}\geq \min_{i=1,\ldots,C(H)}\frac{|S_i|-1}{e(G[S_i])}\geq \min_{S\in Co(G)}\frac{|S|-1}{e(G[S])},$$ where the second inequality holds from Lemma \ref{lemma1}. Thus we can obtain that $$\inf_{\omega\in W} t_\omega(G)\geq\min_{S\in Co(G)}\frac{|S|-1}{e(G[S])}.$$
 Assume $\widetilde{S}\in\arg\min_{S\in Co(G)}\frac{|S|-1}{e(G[S])}.$ Let
 $$\omega'(e)=\left\{
   \begin{array}{ll}
   0 & \mbox{if $e\in E(G)\setminus E(G[\widetilde{S}])$},\\
   1 &  \mbox{if $e\in E(G[\widetilde{S}])$}.
   \end{array}\right.$$
Then
 $$\inf_{\omega\in W} t_\omega(G)  \le  \frac{\omega'(T_{\omega'})}{\omega'(G)} =\frac{|\widetilde{S}|-1}{e(G[\widetilde{S}])} =\min_{S\in Co(G)}\frac{|S|-1}{e(G[S])}.$$
Thus the result holds.  
\qed

\noindent{\bf Remark} By the proof of Theorem \ref{inthm}, there exist $\omega^{\ast} \in \arg\min_{\omega \in W^{0-1}}t_{\omega}(G)$, and $S\in Co(G)$ such that $t_{\omega^{\ast}}(G)=\frac{|S|-1}{e(G[S])}$ with  $\omega^{\ast}(e)=1$ if and only if $e\in E(G[S])$.

In the following subsections, we will characterize the edge-weight function which minimizes the weighting ratio of any chordal graph.
\subsection{The Edge-Weight Function to Minimize Weighting Ratio of a Complete Graph}
We begin with a special chordal graph, the complete graph.

\begin{corollary}\label{complete}
    If $G=K_\kappa$, then $\inf_{\omega \in W} \frac{\omega(T_\omega)}{\omega(G)}=\frac{2}{\kappa}$. Also if $\omega^{\ast} \in \arg\inf_{\omega \in W} \frac{\omega(T_\omega)}{\omega(K_\kappa)}$, then $\omega^{\ast}(e)=1$ for any $e\in E(K_\kappa)$.
\end{corollary}
{\it Proof}  
    Since $G$ is a complete graph, then for any subset $S\subseteq V(G)$, $G[S]$ is also complete. Thus $$\frac{|S|-1}{e(G[S])}=\frac{|S|-1}{{|S| \choose 2}}=\frac{2}{|S|}.$$
  By Theorem \ref{inthm},$$\inf_{\omega\in W} \ratio=\min_{S\subseteq V(G)}\frac{2}{|S|}=\frac{2}{\kappa}.\quad \qed$$
\qed
Recall that a graph is chordal if each of its cycles of length at least 4 has a chord.  The chordal graph has nice properties, one of which is shown below and will be used in Subsection \ref{subsec_3.2}.
\begin{proposition}(\label{clique decompostion}\cite{bondy2008graph})
    Let $G$ be a chordal graph with maximal cliques $C_1,C_2,\ldots,C_k$. Then there exists an permutation $\sigma$ of $\{1,2,\ldots,k\}$ such that $C_{\sigma(t)} \cap (\cup_{s<t} C_{\sigma(s)})$ is a clique for any $t=1,\ldots,k$.
\end{proposition}
\subsection{The Edge-Weight Function to Minimize Weighting Ratio of Any Chordal graph}\label{subsec_3.2}
In this subsection, we propose detailed descriptions about the edge-weight function of chordal graph that minimizes the weighting ratio.
Before that, we define a {\bf clique graph} for any connected graph as follows.
\begin{definition} Let $G$ be a connected graph.
    $G_{c}$ is a clique graph of $G$, if the vertex set $V(G_{c})$ is the collection of maximal cliques in $G$, and can be constructed by following process:

    Step 0. Set $G^{\prime}= \emptyset, G_{c}=\emptyset$.

Step 1. Choose any maximal clique $C$, $C\not\subset G^{\prime}$.

Step 2. If $V(C)\cap V(G') =\emptyset$, then set $G^{\prime}=G^{\prime}\cup C$, $V(G_{c})=V(G_{c})\cup \{v_c\}$ where $v_c$ is a vertex with respect to clique $C$, and go to Step 1. Else, assume that $C$ joint $G^{\prime}$ at $S_1,S_2,\ldots,S_k$ which are the cliques of $G^{\prime}$. For any $1\leq i\leq k$, we choose a maximal clique of $G$, say $C_{i}$, such that $S_i \subseteq C_{i}$. Denote the vertex in $V(G_c)$ with respect to clique $C_i$ by $v_{c_i}$.

Step 3. Set $G^{\prime}=G^{\prime}\cup C$, $V(G_{c})=V(G_{c})\cup \{v_c\}$, and $E(G_{c})=E(G_{c})\cup \{v_cv_{c_{i}}~|~1\le i\le k\}$. If $G^{\prime}= G$, output $G_{c}$, else go to Step 1.
\end{definition}

It should be noticed that we can construct different clique graphs of $G$ if the order of clique selection is different. Let $\widetilde{H}(G)$ denote the set of all clique graphs that can be output by above procedure in different orders of clique selection. If a vertex $v$ of a maximal clique $C$ also belongs to other maximal cliques, then $v$ is called in the public part of $C$; otherwise, we call it in the non-public part.
\begin{theorem}\label{clique_result}
    $G$ is a chordal graph if and only if there exists $G_{c}\in \widetilde{H}(G)$ such that $G_{c}$ is a tree. Moreover, if $G$ is a chordal graph, then for any $G_c\in \widetilde{H}(G)$ and every leaf of $G_c$, the corresponding maximal clique have vertices in non-public part.
\end{theorem}
{\it Proof}  
     Assume $G$ is a chordal graph. From Proposition \ref{clique decompostion}, there is a maximal clique set $\{C_1,\ldots,C_k\}$ with a permutation $\sigma$ that meets the requirement. We will show that
     there is  $G_{c}\in \widetilde{H}(G)$ such that $G_{c}$ is a tree by induction on  $k$. If $k=1$, it is obvious. Let $G^{\prime}=\cup_{i<k}C_{\sigma(i)}$. Then $G^{\prime}$ is a chordal graph with $k-1$ maximal clique. By the inductive hypothesis, there exists $G_{c}^{\prime}\in \widetilde{H}(G^{\prime})$ which is a tree. Since $C_{\sigma(k)}\cap(\cup_{i\leq k-1}C_{\sigma(i)})$ is a clique, one of the clique graph of $G$ can be constructed by adding an one-degree vertex $v_{C_{\sigma(k)}}$ to $G_{c}^{\prime}$. Thus the resulting graph is still a tree.

    Let $G_{c}\in \widetilde{H}(G)$ be a tree with $k$ vertices $v_1,\ldots,v_k$, and the corresponding maximal cliques in $G$ are denoted as $C_1,\ldots,C_k$. We will show that   $G$ is a chordal graph by induction on  $k$.
    We may assume that $v_k$ is a leaf of $G_{c}$. Then $G_{c}-v_k$ is a tree and can be constructed from $G^{\prime}=\cup_{i\leq {k-1}}C_i$. By induction, $G^{\prime}$ is a chordal graph. Then $G$ is also a chordal graph by adding clique $C_k$ to $G^\prime$. By the structure of $G_c$, $C_k$ have vertices in non-public part. 
\qed

The following is the main result of this subsection.
\begin{theorem}\label{chordal ratio}
    Let $G$ be a chordal graph with maximal cliques $C_1,C_2,\cdots, C_k$ satisfying the condition that $C_i \cap (\cup_{j<i}C_j)$ is a clique, and $s$ be the size of maximum clique in $G$. Then there exist a $\omega^{\ast} \in \arg\min_{\omega \in W^{0-1}}t_{\omega}(G)$, and a subgraph $F$ which is induced by $\{e\in E(G)~|~\omega^\ast(e)=1\}$ such that $F$ satisfies the following properties:

    (1) $F$ is a connected chordal graph, and for every tree $F_c\in\widetilde{H}(F)$, all cliques in $F$ corresponding to the leaves of $F_c$ are maximal cliques of $G$;

    (2) $d_F(v)\geq\frac{s}{2}$ for any $v\in V(F)$;

    (3) $\frac{1}{s-1}\leq t_{\omega^\ast}(G)\leq \frac{2}{s}$.
\end{theorem}
{\it Proof}  
   (1)  By induction on  $k$. It is trivial when $k=1$ by Corollary \ref{complete}.  By Remark, there is a $\omega^{\ast} \in \arg\min_{\omega \in W^{0-1}}t_{\omega}(G)$, and a connected chordal subgraph $F$  of $G$ such that $t_{\omega^{\ast}}(G)=\frac{v(F)-1}{e(F)}$ and $\omega^{\ast}(e)=1$ if and only if $e\in E(F)$. Let $S=V(F)$, $D_1,D_2,\ldots,D_r$ be the maximal cliques of $F$ with $d_i :=|D_i|$ and $D_i \subseteq C_i$ ($1\le i\le r$). Then $r\le k$. Assume $\{D_{i_1},\ldots,D_{i_s}\}$ is the set of maximal cliques which corresponds to the leaves of $F_c$, where $F_c\in\widetilde{H}(F)$. Denote $G^{\prime}=\cup _{i=1}^r C_i$. By Theorem \ref{clique_result}, there are vertices in non-public part for any $D_{i_j}$, where $1\le j\le s$. If $r\leq k-1$, by induction hypothesis, there exist a $\omega' \in \arg\min_{\omega \in W^{0-1}}t_{\omega}(G')$ and a connected  chordal subgraph $H$  of $G'$ such that $H$ satisfies property (1). Note that $F$ is a connected of $G'$.  By Theorem \ref{inthm}, we have
     $$t_{\omega'}(G')=\frac{v(H)-1}{e(H)} = \min_{S\in Co(G^{\prime})}\frac{|S|-1}{e(G^\prime[S])}\leq \frac{v(F)-1}{e(F)}=t_{\omega^{\ast}}(G).$$
    Denote $\omega''(e)=\omega'(e)$ if $e\in E(H)$ and $\omega''(e)=0$ if $e\in E(G)\setminus E(H)$. Then $t_{\omega''}(G)=t_{\omega'}(G')\le t_{\omega^{\ast}}(G)$ which implies $ \omega'' \in \arg\min_{\omega \in W^{0-1}}t_{\omega}(G)$ and $H$ is the subgraph we derive.

      So we may assume that $r=k$ and $D_i \cap(\cup_{j<i}D_j)$ is a clique. Let $P_i = D_i\cap(\cup_{j<i}D_j)$ and $p_i=|P_i|$.    Then we can obtain that
    \begin{equation*}
        \ratioohat=\frac{v(F)-1}{e(F)}=\frac{d_1+\cdots+d_k-p_2-\cdots-p_k-1}{{d_1\choose 2}+\cdots+{d_k\choose 2}-{p_2\choose 2}-\cdots-{p_k\choose 2}}.
    \end{equation*}

   \noindent {\bf Claim 1}   $\frac{1}{d_{i_j}}<\frac{1}{d_{i_j}-1}\leq \ratioohat$ for all $1\le j\le s$.

  \noindent{\bf Proof of Claim 1}  Suppose  $\frac{1}{d_{i_1}-1}> \ratioohat$. Let $v$ be a vertex in non-public part of $D_{i_1}$. Set $S'=S\setminus\{v\}$ and $\omega'(e)=1$ if and only if $e\in E(G[S'])$. Then   we have
  $$t_{\omega'}(G)=\frac{\omega^\ast(T_{\omega^\ast})-1}{\omega^\ast(G)-(d_{i_1}-1)}< \ratioohat=t_{\omega^{\ast}}(G),$$ a contradiction.

Suppose there is $j$, say $j=1$, such that $D_{i_1}$ is not a  maximal clique of $G$. Then there is $u\in V(C_{i_1})\setminus V(D_{i_1})$. We first consider the case that $u$ is a vertex in non-public part of $C_{i_1}$. Set $S'=S\cup\{u\}$ and $\omega'(e)=1$ if and only if $e\in E(G[S'])$. By Claim 1,   we have
  $$t_{\omega'}(G)=\frac{\omega^\ast(T_{\omega^\ast})+1}{\omega^\ast(G)+d_{i_1}}< \ratioohat=t_{\omega^{\ast}}(G),$$ a contradiction. Now we  consider the case that $u$ is a vertex in the public part of $C_{i_1}$. Assume $u\in \cap_{j\in J}V(C_j)$, where $J\subseteq \{1,\ldots,k\}$. Let $a_j := |D_j|-|P_j|$ for $j\in J$.
  Set $S'=S\cup\{u\}$ and $\omega'(e)=1$ if and only if $e\in E(G[S'])$. Then
  $$t_{\omega'}(G)=\frac{\omega^\ast(T_{\omega^\ast})+1}{\omega^\ast(G)+d_{i_1}+\sum_{j\in J\setminus\{i_1\}}a_i}\le \frac{\omega^\ast(T_{\omega^\ast})+1}{\omega^\ast(G)+d_{i_1}}.$$By Claim 1, we have $t_{\omega'}(G)<\ratioohat=t_{\omega^{\ast}}(G),$ a contradiction.
  
 (2) Suppose there exists a vertices $v\in V(F)$ satisfying $d_F(v)<\frac{s}{2}$.  By Theorem \ref{inthm}, $t_{\omega^{\ast}}(G)\le  \frac{2}{s}$. Then $\frac{1}{d_F(v)}> t_{\omega^\ast}(G)$. Set $S'=S\setminus\{v\}$ and $\omega'(e)=1$ if and only if $e\in E(G[S'])$. By the same argument as the proof of Claim 1, we can derive a contradiction.

    (3) By Theorem \ref{inthm}, $t_{\omega^{\ast}}(G)\le  \frac{2}{s}$. By Claim 1,  $t_{\omega^\ast}(G)=\ratioohat\ge \frac{1}{d_{i_j}-1}$, where $d_{i_j}=|D_{i_j}$ and $\{D_{i_1},\ldots,D_{i_s}\}$ is the set of maximal cliques which corresponds to the leaves of $F_c$. By (1), we have $t_{\omega^\ast}(G) \ge \frac{1}{s-1}$. 
    \qed
As an application of Theorem \ref{chordal ratio}, we obtain the following conclusion about the maximum chordal subgraph of a graph. Let $t_\kappa(n)$ be the maximum number of edges in the $(\kappa+1)$-clique free graph with $n$ vertices. From Tur\'{a}n Theorem, we have that $t_\kappa(n)=\frac{\kappa-1}{2\kappa}n^2+O(1)$. In order to prove our result, we need the following theorem.
\begin{theorem}(\cite{gishboliner2023maximal}) \label{S} Let $\kappa, n \ge 1$, $G$ a graph of order $n$ and size $m$.  Set $a:= m- t_\kappa(n)$. If $t_\kappa(n)+1\leq m\leq t_{\kappa+1}(n)$, then $G$ contains a chordal subgraph with at least
$$(\kappa-\frac{1}{\kappa})n+\sqrt{\frac{2(\kappa+1)a}{\kappa}}-{\kappa+1\choose 2}-O(\sqrt{n})$$edges.
\end{theorem}

\begin{theorem}
    Let $G$ be a graph with $n$ vertices and $m$ edges such that $t_\kappa(n)+1\leq m\leq t_{\kappa+1}(n)$, $\kappa\geq 2$. Suppose that $F$ is the maximum chordal subgraph of $G$, when $n$ is large enough, there must be a $(\kappa+1)$-clique in $F$.
\end{theorem}
{\it Proof}   Let $s$ be the size of maximum clique in $F$.
    Suppose that $F$ has no $(\kappa+1)$-clique. Then $s\le \kappa$. Let $\omega^{\ast} \in \arg\min_{\omega \in W^{0-1}}\ratioF$. By Theorem \ref{chordal ratio} (3), $t_{\omega^\ast}(F)\geq {1}/{(s-1)}\ge \frac{1}{\kappa-1}$. By Theorems \ref{inthm} and \ref{chordal ratio} (1), we have    \begin{equation*}
         t_{\omega^\ast}(F)\leq \frac{v(F)-1}{e(F)}\leq\frac{n-1}{e(F)}.
    \end{equation*}
    Thus we obtain that
    \begin{equation*}
        e(F)\leq (\kappa-1)(n-1).
    \end{equation*}
    By Theorem \ref{S},
    \begin{equation*}
    e(F)\geq(\kappa-\frac{1}{\kappa})n+\sqrt{\frac{2(\kappa+1)a}{\kappa}}-{\kappa+1\choose 2}-O(\sqrt{n}),
    \end{equation*}
    where $a= m- t_\kappa(n)$. When $n$ is large enough,
    \begin{equation*}
        (\kappa-\frac{1}{\kappa})n+\sqrt{\frac{2(\kappa+1)a}{\kappa}}-{\kappa+1\choose 2}-O(\sqrt{n})> (\kappa-1)(n-1),
    \end{equation*}
    a contradiction.  
\qed
\section{Conclusions}
In this paper, for a given connected graph $G$, we propose an optimal problem,  
whose goal is to find a non-negative edge-weight function on the edges that minimizes the ratio of the weight of its maximum spanning tree to the total weight. Further more we transform this problem into searching the vertex subset $S$ such that $G[S]$ is connected with smallest ratio $\frac{|S|-1}{e(G[S])}$, and prove that this problem is solvable in polynomial time. 


This result can be used in the study of extreme value problems in graph theory. A graph $G$ is 1-balanced if the maximum among $\frac{e(H)}{v(H)-C(H)}$ taken over all subgraphs of $G$, is attained when $H=G$. The results of this paper inspire us to explore the structure about the maximum 1-balanced subgraph of a class of graphs with special properties, which might be useful for further research.

\section*{Acknowledgement} This work is partially supported by the National Natural Science Foundation of China (Grant 12171272).


\begin{thebibliography}{EGOZ89}

\bibitem[AMS91]{DBLP:journals/comgeo/AgarwalMS91}
Pankaj~K. Agarwal, Jir{\'{\i}} Matousek, and Subhash Suri.
\newblock Farthest neighbors, maximum spanning trees and related problems in
  higher dimensions.
\newblock {\em Computational Geometry}, 1:189--201, 1991.

\bibitem[Bal86]{balas1986fast}
Egon Balas.
\newblock A fast algorithm for finding an edge-maximal subgraph with a
  tr-formative coloring.
\newblock {\em Discrete Applied Mathematics}, 15(2-3):123--134, 1986.

\bibitem[BCH15]{bhowmick2015new}
Sanjukta Bhowmick, Tzu-Yi Chen, and Mahantesh Halappanavar.
\newblock A new augmentation based algorithm for extracting maximal chordal
  subgraphs.
\newblock {\em Journal of Parallel and Distributed Computing}, 76:132--144,
  2015.

\bibitem[BM08]{bondy2008graph}
John~Adrian Bondy and Uppaluri Siva~Ramachandra Murty.
\newblock {\em Graph theory}.
\newblock Springer Publishing Company, Incorporated, 2008.

\bibitem[BX24]{DBLP:journals/tcs/BaiX24}
Tian Bai and Mingyu Xiao.
\newblock Exact algorithms for restricted subset feedback vertex set in chordal
  and split graphs.
\newblock {\em Theoretical Computer Science}, 984:114326, 2024.

\bibitem[Cha00]{DBLP:journals/jacm/Chazelle00a}
Bernard Chazelle.
\newblock A minimum spanning tree algorithm with inverse-ackermann type
  complexity.
\newblock {\em Journal of the {ACM}}, 47(6):1028--1047, 2000.

\bibitem[EGOZ89]{erdos1989size}
Paul Erd{\"o}s, Andr{\'a}s Gy{\'a}rf{\'a}s, Edward~T. Ordman, and Yechezkel
  Zalcstein.
\newblock The size of chordal, interval and threshold subgraphs.
\newblock {\em Combinatorica}, 9(3):245--253, 1989.

\bibitem[EK11]{enright2011application}
Jessica~A. Enright and Grzegorz Kondrak.
\newblock The application of chordal graphs to inferring phylogenetic trees of
  languages.
\newblock In {\em Proceedings of the 15th International Joint Conference on
  Natural Language Processing}, pages 545--552, 2011.

\bibitem[Elk06]{DBLP:journals/jcss/Elkin06}
Michael Elkin.
\newblock A faster distributed protocol for constructing a minimum spanning
  tree.
\newblock {\em Journal of Computer and System Sciences}, 72(8):1282--1308,
  2006.

\bibitem[FI03]{fleischer2003push}
Lisa Fleischer and Satoru Iwata.
\newblock A push-relabel framework for submodular function minimization and
  applications to parametric optimization.
\newblock {\em Discrete Applied Mathematics}, 131(2):311--322, 2003.

\bibitem[Gav87]{DBLP:journals/jal/Gavril87}
Fanica Gavril.
\newblock Generating the maximum spanning trees of a weighted graph.
\newblock {\em Journal of Algorithms}, 8(4):592--597, 1987.

\bibitem[GGST86]{DBLP:journals/combinatorica/GabowGST86}
Harold~N. Gabow, Zvi Galil, Thomas~H. Spencer, and Robert~Endre Tarjan.
\newblock Efficient algorithms for finding minimum spanning trees in undirected
  and directed graphs.
\newblock {\em Combinatorica}, 6(2):109--122, 1986.

\bibitem[GMM97]{DBLP:journals/tcs/GalbiatiMM97}
Giulia Galbiati, Angelo Morzenti, and Francesco Maffioli.
\newblock On the approximability of some maximum spanning tree problems.
\newblock {\em Theoretical Computer Science}, 181(1):107--118, 1997.

\bibitem[GS23]{gishboliner2023maximal}
Lior Gishboliner and Benny Sudakov.
\newblock Maximal chordal subgraphs.
\newblock {\em Combinatorics, Probability and Computing}, 32(5):724–741,
  2023.

\bibitem[HK11]{DBLP:journals/jct/HurlbertK11}
Glenn Hurlbert and Vikram Kamat.
\newblock Erd\"{o}s-{K}o-{R}ado theorems for chordal graphs and trees.
\newblock {\em Journal of Combinatorial Theory, Series {A}}, 118(3):829--841,
  2011.

\bibitem[KIM81]{DBLP:journals/siamcomp/KatohIM81}
Naoki Katoh, Toshihide Ibaraki, and Hisashi Mine.
\newblock An algorithm for finding {K} minimum spanning trees.
\newblock {\em {SIAM} Journal on Computing}, 10(2):247--255, 1981.

\bibitem[KKT95]{karger1995randomized}
David~R Karger, Philip~N Klein, and Robert~E Tarjan.
\newblock A randomized linear-time algorithm to find minimum spanning trees.
\newblock {\em Journal of the ACM (JACM)}, 42(2):321--328, 1995.

\bibitem[Kru56]{kruskal1956shortest}
Joseph~B Kruskal.
\newblock On the shortest spanning subtree of a graph and the traveling
  salesman problem.
\newblock {\em Proceedings of the American Mathematical society}, 7(1):48--50,
  1956.

\bibitem[LFJZ18]{DBLP:journals/tcs/LiFJZ18}
Xingfu Li, Haodi Feng, Haitao Jiang, and Binhai Zhu.
\newblock Solving the maximum internal spanning tree problem on interval graphs
  in polynomial time.
\newblock {\em Theoretical Computer Science}, 734:32--37, 2018.

\bibitem[MRG10]{DBLP:conf/kdd/MarchRG10}
William~B. March, Parikshit Ram, and Alexander~G. Gray.
\newblock Fast euclidean minimum spanning tree: algorithm, analysis, and
  applications.
\newblock In {\em Proceedings of the 16th {ACM} {SIGKDD} International
  Conference on Knowledge Discovery and Data Mining}, pages 603--612, 2010.

\bibitem[Pri57]{prim1957shortest}
Robert~Clay Prim.
\newblock Shortest connection networks and some generalizations.
\newblock {\em The Bell System Technical Journal}, 36(6):1389--1401, 1957.

\bibitem[PT22]{DBLP:conf/iwoca/PapadopoulosT22}
Charis Papadopoulos and Spyridon Tzimas.
\newblock Computing a minimum subset feedback vertex set on chordal graphs
  parameterized by leafage.
\newblock In {\em Proceedings of the 33rd International Workshop on
  Combinatorial Algorithms}, pages 466--479, 2022.

\end{thebibliography}
\end{document}